\magnification=\magstep1

\newcount\sec \sec=0
\input Ref.macros
\input math.macros
\input labelfig.tex
\input epsf
\forwardreferencetrue
\citationgenerationtrue
\initialeqmacro


\title{Approximating Cayley diagrams versus Cayley graphs
}

\author{\'Ad\'am Tim\'ar}
\bigskip

\abstract{We construct a sequence of finite graphs that weakly converge to
a Cayley graph, but there is no labelling of the edges that would converge
to the corresponding Cayley diagram. A similar construction is used to 
give graph
sequences that converge to the same limit, and such that a Hamiltonian cycle 
in one of them has a limit that is not approximable by any subgraph of the 
other. We give an
example where this holds, but convergence is
meant in a stronger sense. This is related to whether having a
Hamiltonian cycle is a testable graph property.}


By a {\it diagram} we will mean a graph with edges oriented and labelled by
elements of some given set (of ``colors").
By a {\it Cayley diagram} we mean a Cayley graph when we do not forget that edges are oriented and labelled by
elements of the generating set.
For simplicity, when an edge is oriented both ways with the same label
(that is, when the labelling generator has degree 2), we will
represent it by an unoriented edge in the Cayley diagram, and
sometimes refer to it as a 2-cycle. A {\it rooted graph (diagram)} is a graph (diagram) with a distinguised vertex called the root; a rooted isomorphism between rooted graphs $G$ and $H$ is an isomorphism that maps root into root. A rooted labelled-isomorphism between rooted diagrams $G$ and $H$ is a rooted isomorphism that preserves orientations and labels of the edges.

Let ${\cal G}$ be the set of rooted isomorphism classes of countable connected rooted graphs. Let $
\tilde {\cal G}$ be the set of rooted isomorphism classes of connected rooted diagrams. We can introduce a metric on ${\cal G}$ by saying that the distace between $G,H\in{\cal G}$ is $2^{-r}$ if $r$ is the largest integer such that the $r$-neighborhood of the root in $G$ is rooted isomorphic to the $r$-neighborhood of the root in $H$. We can define distance on $\tilde {\cal G}$ similarly. It is easy to check that the generated topology makes ${\cal G}$ ($\tilde {\cal G}$) into a complete separable metric space.

Suppose that $G_n$ is a sequence of finite graphs. Then we can define a probability measure $\mu_n$ on ${\cal G}$
by picking a vertex $o$ uniformly at random as the root, and projecting the resulting measure to ${\cal G}$. Now, say that $G_n$ {\it converges} to a probability measure $\mu$ on ${\cal G}$, if the $\mu_n$ weakly converges to $\mu$ (i.e., for each bounded continuous function $f:{\cal G}\to \R$ we have $\int_{\cal G} f d\mu_n \to\int_{\cal G} f d\mu$). This convergence is often called Benjamini-Schramm convergence; see \ref b.AS/ or \ref b.AL/ for more details.
If $G$ is some transitive graph, then we can define a Dirac delta measure $\mu$ on ${\cal G}$ that is supported on the rooted isomorphism class of $(G,o)$, where $o$ is an arbitrary point. If $G_n$ converges to this $\mu$, then we will say that $G_n$ {\t converges to} $G$, or that $G$ {\t is approximated by} $G_n$.
Similarly, if $G$ is quasi-transitive, there is a natural finitely supported probability measure on $\{(G,o_1),\ldots, (G,o_m)\}$, where $\{o_1,\ldots, o_m\}$ in $G$ is a traversal for the orbits of the automorphism group of $G$. The same definitions and terminology apply for diagrams instead of graphs (where ``rooted isomorphism'' is replaced by ``rooted labelled-isomorphism''). 

Less formally, convergence of $G_n$ to a transitive $G$ means that for any $r$, the proportion of vertices $x$ in $G_n$ whose $r$-neighborhood with $x$ as a root is rooted isomorphic to the $r$-neighborhood of $o$ in $G$ tends to 1 as $n\to \infty$. 
It is a central open question whether any {\it unimodular} transitive 
graph 
can be approximated by a sequence of finite graphs. See \ref b.AL/ for the 
definition of unimodularity (which is a necessary condition for the 
existence of such an approximation), and for more details on what we have introduced. Cayley graphs are unimodular, and a finitely generated group is called {\it sofic}, if it has a finitely generated Cayley diagram that is approximable by a sequence of finite diagrams. (There are several equivalent definitions of soficity, see \ref b.ESz/ or \ref b.P/ for history and references.) The interest in whether every group is sofic comes partly from the fact that many conjectures are known to hold for sofic groups. A nice brief survey on the subject is \ref b.P/.

By definition, if a sequence $\tilde G_n$ of finite diagrams converges to a 
Cayley diagram $\tilde G$, then the underlying graphs $G_n$ converge to the 
underlying Cayley graph $G$. It is natural to ask, whether the converse is 
true, or whether the approximability of a Cayley graph by finite graphs 
implies that the group is sofic. The next two questions phrase this in increasing difficulty. The second one seems to have been asked by several people independently. The first one was proposed by Russell Lyons at a workshop in Banff \ref b.ASz/.

\procl q.1
Suppose that a sequence $G_n$ of finite graphs converges to a Cayley graph $G$, and let $\tilde G$ be a Cayley diagram with underlying graph $G$. Is there a sequence of diagrams $\tilde G_n$ such that if we forget about orientations and labels of edges in $\tilde G_n$ we get $G_n$, and such that the sequence $\tilde G_n$ converges to the diagram $\tilde G$?
\endprocl

\procl q.2 
Suppose that a Cayley graph of a group $\Gamma$ is approximable by a sequence of finite graphs. Is $\Gamma$ then sofic?
\endprocl

We give a negative answer to \ref q.1/ in \ref t.main/. This indicates 
that the existence of an approximating sequence for a Cayley graph may not 
help directly in the construction for an approximation of the Cayley 
diagram. In fact, it is reasonable to think that to answer \ref q.2/ might 
be as difficult as the question whether every group is sofic. The difficulty of \ref q.2/ is further illustrated by the fact, as explained to us by G\'abor Elek, that some Burger-Mozes groups are known to have a Cayley graph that is the direct product of two regular trees, even though these groups are simple and not known to be sofic. The product of trees is clearly approximable by a sequence of finite graphs, hence a positive answer to \ref q.2/ would imply that these groups are sofic. See IV.9. in \ref b.dlH/ for more on isometric Cayley graphs and Burger-Mozes groups.

If a sequence of graphs (diagrams) $G_n$ weakly converges to a graph
(diagram) $G$, we will write $G_n\to G$. Given a graph or diagram $G$
and vertex $v\in V(G)$, we denote the $r$-neighborhood of $v$ in $G$
by $B_G (v, r)$.

\procl t.main
There exists a Cayley diagram $\tilde G$ such that the corresponding
Cayley graph $G$ is the weak limit of a sequence $G_n$ of randomly
rooted finite graphs, but there is no sequence of diagrams $\tilde
G_n$ that would weakly converge to $\tilde G$ and such that the graph
underlying $\tilde G_n$ is $G_n$.
\endprocl

\proof
Consider $G=T\times{\bf  C_4}$, where $T$ is the 3-regular tree, $C_4$ is
the cycle of length 4, and in the direct product two edges are
adjacent by definition iff they are equal in one coordinate, and
adjacent in the other.
In other words, we have four copies $T_1$, $T_2$, $T_3$, $T_4$ of the
3-regular tree (that we will also call {\it fibers}), some isomorphisms $\phi_1 : T_1\to T_2$, $\phi_2 :
T_2\to T_3$, $\phi_3 : T_3\to T_4$, $\phi_4 : T_4\to T_1$ such that
$\phi_4^{-1}=\phi_3\circ \phi_2\circ \phi_1$; and $G$ consists of
$T_1\cup T_2\cup T_3\cup T_4\cup K$, where $K$ denotes the set of all
edges of the form $\{v, \phi_i (v)\}$ (in particular, $K$ consists of
cycles of length 4).

Let $\tilde G$ be the following diagram. We consider the Cayley
diagram $\tilde T$ of $\Z * \Z _2=\langle a,b | b^2\rangle$. Make $T_i$ a Cayley diagram
labelled-isomorphic to $\tilde T$, and do it in such a way that the
$\phi_i$ are labelled-isomorphisms. To define labels on elements of
$K$, we will use colors $c$ and $d$. Namely, for each 4-cycle in $K$,
color the edges by $c$ and $d$ alternatingly. Do it in such a way that
if the edge between $v$ and $\phi_i (v)$ has label $c$, then for all
neighbors $w$ of $v$ in $T_i$, the edge $\{w, \phi_i (w)\}$ will have
color $d$; and similarly with $c$ and $d$ interchanged.

We claim that the resulting $\tilde G$ is a Cayley diagram. Consider
$$\langle a,b,c,d| b^2, c^2, d^2, cdcd, ada^{-1}c, aca^{-1}d, bcbd\rangle.$$
To see
that the corresponding Cayley diagram is indeed the diagram $\tilde G$ that we
defined, note that the latter has a cycle space generated by 2- and 4-cycles.
The relators given here, together with some of their conjugates of reduced lengths 4, are exactly the
words read along 2- and 4-cycles on a given vertex. 

Now, let $H_n$ be a sequence of 3-regular graphs with girth tending to
infinity and independence ratio less than $1/2-\epsilon<1/2$. See \ref b.B/ for the construction of such a sequence (with $\epsilon=1/26$). Define $G_n=H_n\times C_4$. Clearly,
$G_n\to G$.
Suppose now, that there is a diagram $\tilde G_n$ with underlying
graph $G_n$ such that $\tilde G_n$ weakly converges to $\tilde G$. Let
$H_n^1$, $H_n^2$, $H_n^3$, $H_n^4$ be the four copies of $\tilde H_n$
in $\tilde G_n$ (we will call them {\it fibers} of $\tilde G_n$). Fix a point $o$ of $G$ in $T_1$.
Say that $x\in \tilde G_n$ is {\it $R$-good}, if
there exists a rooted labelled-isomorphism from $B_{\tilde G_n} (x,R)$
to $B_{\tilde G} (o,R)$. 

For $R\geq 4$, the ball $B_{\tilde G} (o,R)$ has only one
rooted labelled-isomorphism to itself, the identity. This is so because every
rooted labelled-isomorphism has to preserve edges in $T_1\cup T_2\cup T_3 \cup
T_4$, and the
only two rooted labelled-isomorphisms that respect the labels and
orientations on $T_1,T_2,T_3,T_4$ can be the identity and one that
switches $T_2$ and $T_4$. The latter, however, switches edges of labels
$c$ and $d$, hence it is not a rooted labelled-isomorphism. 
As a consequence of the fact just proved, if a graph is rooted labelled-isomorphic to 
$B_{\tilde G} (o,R)$ ($R\geq 4$), then there is a unique isomorphism between them.
For each $R$-good $x$, each
$i$ and 
each $\iota$
rooted labelled-isomorphism from $B_{\tilde G_n} (x,R)$ to $B_{\tilde G} (o,R)$, 
every vertex of $B_{G_n}(x,R)\cap
H_n^i$ is mapped into the same $T_j$ by $\iota$ (that is, if two point are in the same fiber, then they are mapped into the same fiber by the rooted labelled-isomorphism). This is so because preserving labels on the edges means in particular that edges within a fiber (of label $a$ or $b$) are mapped into edges within a fiber (the ones having label $a$ or $b$).
There is at most one such rooted labelled-isomorphism
(since if there were more, that would give a nontrivial rooted
labelled-isomorphism from $B_{\tilde G} (o,R)$ to itself, as observed above). We have
obtained that for every $R$-good $x\in \tilde G_n$ ($R\geq 4$), there
is a unique rooted labelled-isomorphism $\iota_x$ from $B_{\tilde G_n} (x,R)$ to
$B_{\tilde G} (o,R)$, and it maps fibers into fibers (in a bijective way).
\def\Sl{\overleftarrow S_n}
\def\Sr{\overrightarrow S_n}

Since $\iota_x$ preserves fibers and is an isomorphism, it either changes the cyclic
order of $H_n^1,H_n^2,H_n^3,H_n^4$ (meaning $\iota_x (H_n^1\cap
B_{\tilde G_n} (x,R)) =B_G(o,R)\cap T_j, 
\iota_x (H_n^2\cap
B_{\tilde G_n} (x,R)) =B_G(o,R)\cap T_{j-1}, \ldots$), or preserves the
cyclic orientation. Let $\Sr$ be the set of $R$-good points $x$ in $\tilde
G_n$ where $\iota_x$ preserves the cyclic order, and $\Sl$ be the set
of those where it reverses the cyclic order. We claim that if $x$ and $y$ are $R$-good and
adjacent in $\tilde G_n$, then $\iota _x$ and $\iota_y$ give different
orientations. To see this, let the $c$-edges adjacent to $x$ and $y$ in $\tilde G_n$ be $\{x,x'\}$ and $\{y,y'\}$
respectively. Then $\iota_x (x')$ and $\iota_x (y')$ are in different fibers, 
hence $x'$ and $y'$ are in different fibers too. On the other hand
$\iota_x(x')$ and $\iota_y (y')$ are in the same fiber by definition, hence one of $\iota_x$ and $\iota_y$ has to preserve orientation and the other one has to reverse it.

We conclude that $\Sr$ is an
independent set, and also $\Sl$ is an independent set. By the choice
of the $H_n$ we then have
$|\Sr \cap H_n^i|/|H_n^i|\leq {1\over 2}-\epsilon$ for every $i$, and
similarly for the $\Sl$. Hence

$$|\Sr\cup\Sl |/|\tilde G_n|=\sum_{i=1}^4 |\Sr\cap
H_n^i|/4|H_n|+|\Sl\cap H_n^i|/4|H_n|\leq 1-2\epsilon.$$
This is uniform in $n$, contradicting the fact that the proportion of
$R$-good points in $\tilde G_n$ (that is, $\Sr\cup\Sl$) 
tends to 1.
\Qed

G\'abor Elek has asked the following question. A positive answer would show that having a Hamiltonian cycle is a testable graph property. (A property being testable is, vaguely, the following. Given a finite graph $G$, can we decide by sampling a bounded number of balls in it, whether there is a graph $G'$ with the property in question, and such that one can transform $G$ into $G'$ by changing an at most $\epsilon$ proportion of the edges in $G$? See \ref b.L/ for the precise definition) . A result of this type is the one in \ref b.EL/, where it is shown that for a convergent graph sequence the matching ratio (that is, the ratio of the size of a maximal matching and the size of the graph) also has a limit. This implies that the matching ratio is a testable graph parameter, \ref b.E/. See \ref b.L/ for the relevance of parameter testing and its connection to graph sequences.

\procl q.4
Let $G_n$ and $H_n$ be two graph sequences, converging to the same (random) $G$. Suppose that $G_n$ contains a Hamiltonian cycle $C_n$ (whose limit is then a biinfinite path). Is there a subgraph in $H_n$ whose limit is the same?
\endprocl

We construct an example where convergence to the same limit fails only in a stronger sense, namely, that there is no subgraph $D_n$ in $H_n$ such that the pair $(H_n, D_n)$ would converge to the same pair, as $(G_n, C_n)$. 
(For $C_n$ subgraph of $G_n$ on the same vertex set, one can think about the pair $(G_n,C_n)$ as a diagram on $G_n$, simply by coloring edges of $C_n$ to one color and edges outside of $C_n$ to another one. Convergence of the pairs $(G_n,C_n)$ can then be defined as convergence of the respective diagrams.)
Our example will use the one given in \ref t.main/.

\procl t.elek
There are two sequences, $G_n$ and $K_n$, that converge to the same Cayley graph $G$, and such that $K_n$ contains a Hamiltonian cycle $C_n$ such that the pair $(K_n,C_n)$ converges to $(G, F)$, where $F\subset G$ is a unimodular random graph, but $G_n$ does not have any subgraph $D_n$ such that $(
G_n, D_n)$ would converge to $(G,F)$.

\endprocl

\proof
Consider $G=T\times C_4$ as in \ref t.main/.

Let $G_n=H_n\times C_4$ be as in \ref t.main/.
We have seen that the limit of $G_n$ is $G$.


The other sequence, $K_n$ will also be the direct product of a graph $B_n$ and $C_4$, and it will have the property that it contains a Hamiltonian cycle such that every other edge of the Hamiltonian cycle is an edge ``coming from $C_4$ (by which we mean an edge of the form $\{(x,v),(x,w)\}$, $x\in B_n$, $v,w\in C_4$ adjacent).
So, consider a bipartite graph $B_n$, with the following properties. It is 3-regular, it contains a Hamiltonian cycle, it has ``upper set'' $U_n=U$ and ``lower set'' $L_n=L$ both containing $2n+1$ vertices, and the girth tends to infinity as $n\to\infty$. We will construct $K_n$ as follows. First, define a bipartite {\it directed} graph $K_n'$ on vertex set $V_1\cup V_2\cup V_3\cup V_4$ where $V_1=\{x_v^1\,:\,v\in U\}$, $V_3=\{x_v^3\,:\,v\in U\}$, $V_2=\{x_w^2\,:\,w\in L\}$, $V_4=\{x_w^4\,:\,w\in L\}$, and set of directed edges $\{(x_s^i,x_t^{i+1})\,:\, \{s,t\}\in E(B_n)   \}$, where $i+1$ is modulo 4 (and similarly later for such indices, without further mention). That is, for each pair $V_i,V_{i+1}$, we ``copy'' $B_n$ on $V_i\cup V_{i+1}$, ($V_i$ playing the role of $U$ iff $i$ is odd), and orient the edges from $V_i$ towards $V_{i+1}$. In particular, $K_n'$ has $2(4n+2)$ vertices, each having indegree 3 and outdegree 3, and all edges going out of $V_i$ go to $V_{i+1}$. To finish, let $K_n=H$ be a bipartite graph of $4(4n+2)$ vertices, whose vertex set is obtained by doubling every vertex $w$ of $K_n'$ to get the twins $\bar w$, $\hat w$. Let $\bar w$ and $\hat v$ be adjacent in $K_n$ iff there is a (directed) edge from $w$ to $v$ in $K_n'$. Further, connect each pair of twins $\bar w$, $\hat w$ by an edge, and call the edges of this type {\it blue} edges. Finally, if $x_v^1\in V_1$, $x_v^3\in V_3$ (with a $v\in V(U_n)$), then connect $\bar x_v^1$ and $\hat x_v^3$ by an edge, and connect $\hat x_v^1$ and $\bar x_v^3$ by an edge. Similarly, if $x_v^2\in V_2$, $x_v^4\in V_4$ (with a $v\in V(L_n)$), then connect $\bar x_v^2$ and $\hat x_v^4$ by an edge, and connect $\hat x_v^2$ and $\bar x_v^4$ by an edge. Call the edges of these type {\it yellow} edges. Observe that the colored edges of
$K_n$ form cycles of lengths 4, each colored by yellow and blue alternatingly. More precisely, note that $K_n$ is isomorphic to $B_n\times C_4$. To see this, note that each of the four sets $\{\bar x_v^i \,:\, v\in U_i  \}\cup\{\hat x_w^{i+1}  \,:\, w\in L_{i+1}   \}$, $i=1,3$, and $\{\bar x_v^i \,:\, v\in L_i  \}\cup\{\hat x_w^{i+1}  \,:\, w\in U_{i+1}   \}$, $i=2,4$
induces a graph isomorphic to $B_n$. We will refer to these four sets as fibers. Colored edges of $K_n$ then correspond to edges coming from $C_4$ in the direct product. In particular, it is clear that $K_n$ converges to $G$. Consider now the set $S$ of blue edges with one endpoint in a fiber $C^1$ and the other endpoint in fiber $C^2$. In the direct product, the  4-cycles that correspond to neighboring vertices have alternating colorings, hence the endpoints of $S$ in $C^1$ form an independent set (since $B_n$ is bipartite). We will refer to this as the ``independence property''.

We claim that $K_n$ contains a Hamiltonian cycle. To see this, let the vertices in a Hamiltonian cycle of $B_n$ be $v_1,v_2,\ldots, v_{4n+2}$, listed in their order along the cycle. The respective vertices $x_{v_1}^1,x_{v_2}^2,x_{v_3}^3,v_{v_4}^4,x_{v_5}^1,x_{v_6}^2,\ldots,x_{v_{4n+2}}^2, x_{v_1}^3,x_{v_2}^4,x_{v_3}^1,x_{v_4}^2,\ldots, x_{v_{4n+2}}^4$ determine a Hamiltonian directed cycle in $K_n'$. These edges can be projected into $K_n$, and if we add the (blue) edge between each pair $\bar x_v^i$, $\hat x_v^i$, we get a Hamiltonian cycle $C_n$ of $K_n$. Every second edge on $C_n$ is blue. 

Now, let $\Omega$ be the set of edges of $G$ not in the fibers (that is, edges coming from $C_4$). We have seen that local isomorphisms from $K_n$ to $G$ map colored edges to edges in $\Omega$. Hence the limit of $C_n$ in $G$ is a biinfinite path $F$ that has every other edge in $\Omega$. Fibers are also preserved, thus by the ``independence property'' we obtain for the set of edges of $F\cap \Omega$ with one endpoint in a fiber $C^1$ and the other in a fiber $C^2$, that the set of their endpoints in $C^1$ is independent.

Suppose now that there is a subgraph $D_n\subset G_n$ such that $(G_n, D_n)$ would converge to $(G,F)$. We proceed similarly as in the proof of \ref t.main/. Fix $o\in V(G)$, and let $X$ be the set of $R$-good points $x$ such that the (unique) local isomorphism from $B_{G_n}(x,R)$ to $B_G(o,R)$ does not change the (previously fixed) orientation of the fibers. By the same argument as in the last two paragraphs of the proof of \ref t.main/, $X$ is an independent set. Furthermore, its density is larger than $(1-\epsilon)/2$ if $n$ is large enough, since $G_n\to G$. This contradicts the assumption on the size of the largest independence set in $G_n$.


\Qed
\bigskip
\SetLabels
   (.17*.74) {$v$} \\ 
   (.174*.4) {$w$} \\
  (.15*.85) {$U$}\\
  (.15*.27) {$L$}\\
  (.35*.72)  {$x_v^1$}\\
  (.35*.4)  {$x_w^2$}\\
  (.45*.72)  {$x_v^3$}\\
  (.45*.4)  {$x_w^4$}\\
  (.33*.85) {$V_1$}\\
  (.33*.27) {$V_2$}\\
  (.43*.85) {$V_3$}\\
  (.43*.27) {$V_4$}\\
  (.87*.71) {$\bar x_v^3$}\\
  (.87*.4) {$\hat x_w^4$}\\
  (.5*.07) {{\sl The scheme of constructing $K_n'$ and 
$K_n$ from $B_n$.}}\\
\endSetLabels
\AffixLabels{\epsfysize=5.3cm
\epsfbox{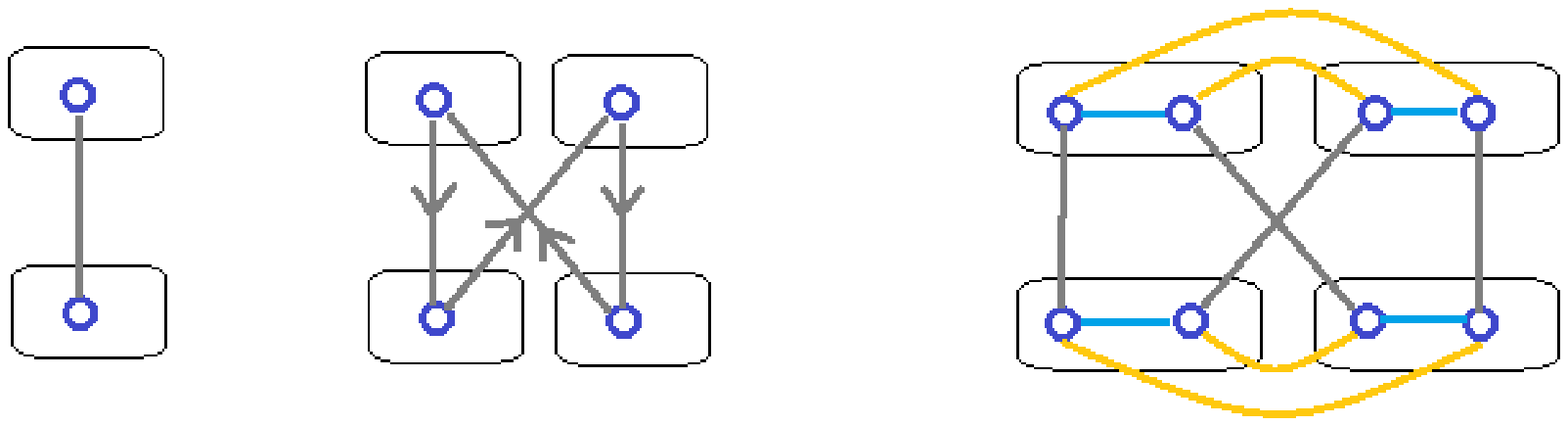}
}


\medbreak
\noindent {\bf Acknowledgements.}\enspace
I am very grateful to Russell Lyons and G\'abor Elek for motivating 
questions and helpful discussions. I also thank G\'abor Lippner for his 
comments on the manuscript.

Research supported by Sinergia grant CRSI22-130435 of the Swiss
National Science foundation.

\startbib

\bibitem[ASz]{ASz} Ab\'ert. M. \and Szegedy, B. (2009) Residually finite groups, graph limits and dynamics, Banff reports, available at 
http://www.birs.ca/workshops/2009/09frg147/report09frg147.pdf

\bibitem[AL]{AL} Aldous, D. \and Lyons, R. (2007) Processes on unimodular 
random
networks {\it Electr. Commun. Prob.} {\bf 12}, 1454-1508.

\bibitem[AS]{AS} Aldous, D. \and J.M. Steele (2004) The Objective Method:
Probabilistic Combinatorial Optimization
and Local Weak Convergence.
In Kesten, H., editor, {\it Probability on
Discrete Structures}, volume 110 of {\it Encyclopaedia Math. Sci.}, 1-72. Springer,
Berlin. Probability Theory, 1.

\bibitem[B]{B} Bollob\'as, B. (1981) The independence ratio of regular graphs {\it Proc. Amer. Math. Soc.}
{\bf 83}, 433-436.

\bibitem[dlH]{dlH} P. de la Harpe. {\it Topics in Geometric Group Theory}. 
Chicago Lectures in Mathematics Series, 2000.

\bibitem[E]{E} Elek, G. (2010) Parameter testing in bounded degree graphs of subexponential growth {\it Random Structures and Algorithms} {\bf 37}, 248-270. 

\bibitem[EL]{EL} Elek, G. \and Lippner G. (2010) Borel oracles. An analytic approach to constant time algorithms {\it Proc. Amer. Math. Soc.} {\bf 138}, 2939-2947.

\bibitem[ESz]{ESz} Elek, G. \and Szab\'o, E. (2005) Hyperlinearity, essentially free actions and L2-invariants. The sofic
property {\it Math. Ann.} {\bf 332}, no. 2, 421?441.

\bibitem[L]{L} Lov\'asz, L. (2009) Very large graphs. In {\it Current Developments in Mathematics} Vol. 2008, 67-128.

\bibitem[P]{P} Pestov, V.G. (2008) Hyperlinear and sofic groups: a brief guide {\it The Bulletin of Symbolic Logic} {\bf 14}, no. 4, 449-480.

\endbib

\bibfile{\jobname}
\def\noop#1{\relax}
\input \jobname.bbl

\filbreak
\begingroup
\eightpoint\sc
\parindent=0pt\baselineskip=10pt
Fakult\"at f\"ur Mathematik, Universit\"at Wien
Nordbergstrasse 15, 1090 Wien
\endgroup

\bye